\documentclass[aps]{revtex4}

\usepackage{epsfig}
\usepackage{amsmath,amssymb,color}
\usepackage{footnote}
\usepackage[english]{babel}

\parskip=\medskipamount

\newcommand{\eq}[1]{(\ref{#1})}

\newcommand{\be}{\begin{equation}}
\newcommand{\ee}{\end{equation}}

\newcommand\disp{\displaystyle}

\begin{document}

\title{Brownian flights over a circle}

\author{Alexander Vladimirov$^{1}$, Senya Shlosman$^{1,2,3}$, and Sergei
Nechaev$^{4,5}$}

\affiliation{$^{1}$Institute of Information Transmission Problems RAS, 127051 Moscow,
Russia \\ $^{2}$Skolkovo Institute of Science and Technology, 143005 Skolkovo, Russia \\
$^{3}$Aix-Marseille University, Universite of Toulon, CNRS, CPT UMR 7332, 13288, Marseille,
France \\ $^{4}$Interdisciplinary Scientific Center Poncelet, CNRS UMI 2615, 119002 Moscow, Russia \\$^{5}$P.N. Lebedev Physical Institute RAS, 119991 Moscow, Russia}

\begin{abstract}

The stationary radial distribution, $P(\rho)$, of the random walk with the diffusion coefficient $D$, which winds with the tangential velocity $V$ around the impenetrable disc of radius $R$ for $R\gg 1$ converges to the distribution involving Airy function. Typical trajectories are localized in the circular strip $[R, R+ \delta R^{1/3}]$, where $\delta$ is the constant which depends on the parameters $D$ and $V$ and is independent on $R$.

\end{abstract}

\maketitle

\section{Introduction}

We study the behavior of the Brownian motion $B$ in $\mathbb{R}^{2}$, evading a certain region. The questions of similar nature have been considered in the literature. A notable example is the paper \cite{FS} by P. Ferrari and H. Spohn, in which they discussed the limiting behavior of the directed Brownian bridge in (1+1) dimensions constrained to stay above the semicircular obstacle $D_{R}$ of the radius $R.$ It has been found in \cite{FS} that the constrained Brownian trajectory in the limit $R\rightarrow\infty$ converges to a certain diffusion process, which stays at a typical distance $\sim R^{1/3}$ above $D_{R}$ and after rescaling by $R^{1/3}$ has a stationary distribution expressed in terms of the Airy function. In the later paper (\cite{ISV}) such a diffusion process (living in the halfspace $\mathbb{R}_{+}^{2}$) was named "Ferrari-Spohn diffusion".

Recently, in \cite{NPSVV} we studied a related problem, in which the Brownian bridge is replaced by the true 2D Brownian motion, avoiding a disk (or a triangle) in $\mathbb{R}^{2}.$ The following question was the subject of the paper \cite{NPSVV}: could a two-dimensional random path pushed by some constraints to an improbable "large deviation regime", possess extreme statistics with one-dimensional Kardar-Parisi-Zhang (KPZ) fluctuations? It has been shown that the answer is positive, though non-universal, since fluctuations depend on the underlying geometry. In \cite{NPSVV} two examples of 2D systems have been considered in details: the semicircle and the triangle, for which imposed external constraints forced the underlying stationary stochastic process to stay in an atypical regime with anomalous statistics. In the paper \cite{HP} the question of the behavior of a planar Brownian loop capturing a large area $A$ is studied. In particular, it has been established that a typical loop with a fixed area avoids a disc of area $\sim A$, as $A\to\infty$.

In the present paper we study a behavior of the Brownian motion $B$ in $\mathbb{R}^{2}$ which cannot enter a disc of radius $R$, moving around it with a given linear (tangential) velocity. We look at that problem from the viewpoint of the grand canonical ensemble. First of all, we find the drift field $\mathbf{v}$ (\textit{the tilt}) under which the (unconstrained) Brownian particle behaves as an atypical Brownian particle winding around the disc $D_{R}$. Secondly, we determine the stationary distribution of such a particle. We find that the corresponding distribution, scaled by $R^{1/3},$ is expressed via the square of the Airy function (as $R\to\infty$) -- see Eq.(\ref{eq:23}), which is the main result of our paper.

\section{Statement of the problem}

\begin{figure}[th]
\epsfig{file=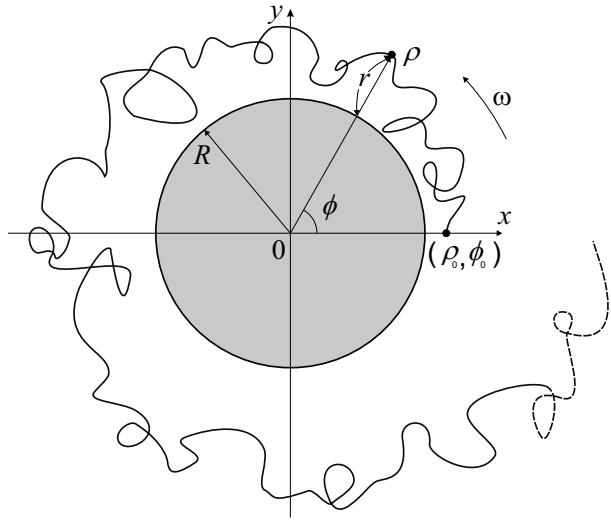,width=8cm}\caption{Schematic image of a particular realization of a two-dimensional random walk with fixed average angular velocity $\omega$ above the circle. The path starts at the point $(\rho_{0},\phi_{0})$. We are interested in finding the stationary time- and
angle-independent distribution $P(\rho)$ for ensemble of such trajectories.}
\label{fig:01}
\end{figure}

We study large deviation problem for a biased two-dimensional Brownian motion $B$ on a plane in a radial geometry. We are looking for the rotation-invariant drift field $\mathbf{v}$, which imposes the following constraints: i) the Brownian motion always stays outside the disk $D_{R}$ of radius $R$, and ii) the Brownian motion winds around the disc with the mean angular velocity $\omega= \frac{V}{2\pi R}$ where the mean linear velocity $V$ is fixed. The requested drift field $\mathbf{v}$ should produce the same behavior as that of the Brownian particle constrained to stay away from $D_{R}$ and making a full turn around $D_{R}$ during typical time $\frac{2\pi R}{V}$. In other words, we are looking for the angular tilt which would produce the constrained behavior as its typical one. Schematically, the considered system is shown on Fig. \ref{fig:01}.

Our goal is to find the stationary distribution $P(\rho,\phi)$ of the process described above, where by $(\rho,\phi)$ we denote polar coordinates on the plane. Clearly, $P(\rho,\phi)= P\left(  \rho\right)$, and the rotation-invariant drift field $\mathbf{v}$ is specified by its projection to one ray.

The time-dependent distribution $P(\rho,\phi,t)$ of the Brownian motion in the rotation-invariant drift field $\mathbf{v}$ satisfies the nonstationary Fokker-Planck (FP) equation, which in polar coordinates can be written as follows:
\be
\frac{\partial P}{\partial t}=D\left[\frac{1}{\rho}\frac{\partial} {\partial\rho} \left(\rho\frac{\partial P}{\partial\rho}\right)+\frac{1}{\rho^{2}} \frac{\partial^{2}P} {\partial\phi^{2}}\right]-\frac{1}{\rho}\frac{\partial}{\partial\rho}\left(\rho v_{\rho}P\right)-\frac{1}{\rho}\frac{\partial}{\partial\phi} \left(v_{\phi}P\right)
\label{eq:01}
\ee
where $v_{\phi},v_{\rho}$ are polar coordinates of the vector $\mathbf{v}(\phi,\rho)$ and $D$ is the diffusion coefficient.

Suppose the drift field $\mathbf{v}$ to be such that the corresponding diffusion $P(\rho,\phi,t)$ stays away from the disc of radius $R$, is ergodic with a stationary distribution $P(\rho,\phi)= P\left(\rho\right)$, and the prescribed mean linear velocity is $V$. Denote by $\Omega(R,V)$ the set of all corresponding drift fields $\mathbf{v}$. According to the Schilder theorem \cite{DZ,R}, the large deviation rate function for the Brownian motion is the square of the velocity. Hence the drift field $\mathbf{v}\in\Omega(R,V)$ we are looking for, is the one which minimizes the action $S$ (the "entropy" of ensemble of Brownian trajectories)
\be
S(\mathbf{v})=\int_{R}^{\infty}\mathbf{v}^{2}(\rho,\phi)P_{\mathbf{v}}(\rho,\phi)\, \rho d\rho\,d\phi
\label{eq:00}
\ee
over $\Omega(R,V)$.

First of all, we obtain a closed expression for the action $S$ in terms of $P(\rho)$ only, expressing $v_{\rho}$ and $v_{\phi}$ as functionals of $P(\rho)$ by solving separately two auxiliary problems. Then we study the limiting behavior of the stationary distribution $P(\rho,\phi)$ as $R$ tends to the infinity, and we see the emergence of the universal distribution expressed terms of the Airy function in that limit.

\subsection{Determination of radial velocity, $v_{\rho}(\rho)$}

Let us point out that if $\mathbf{v}$ is rotation-invariant and $P$ is a stationary rotation-invariant solution of the equation (\ref{eq:01}), then the radial velocity $v_{\rho}(\rho)$ is defined by $P$. Indeed, the distribution $P$ and velocity $v_{\phi}$ do not depend on $t$ and $\phi$, so the Fokker-Planck (FP) equation (\ref{eq:01}) boils down to
\be
D\frac{\partial}{\partial\rho}\left(\rho\frac{\partial P(\rho)}{\partial\rho}\right)- \frac{\partial}{\partial\rho}\big(\rho v_{\rho}P(\rho)\big)=0
\label{eq:03}
\ee

Integrating the FP equation (\ref{eq:03}), we get
\be
D\frac{\partial P(\rho)}{\partial\rho}=v_{\rho}P(\rho)+\frac{C_{0}}{\rho},
\label{eq:04}
\ee
where $C_{0}$ is the integration constant, yet undetermined. Recall that $P$ is a stationary rotation-invariant distribution for which the difference $D\frac{\partial P(\rho)}{\partial\rho}-v_{\rho}P(\rho)$ vanishes since it equals to the mass transfer through the circle of radius $\rho$. Thus, the mass transfer is zero in the stationary case which implies $C_{0}=0$, and for the radial velocity one gets the expression
\be
v_{\rho}(\rho)=\frac{D}{P(\rho)}\frac{\partial P(\rho)}{\partial\rho}
\label{eq:05}
\ee

\subsection{Determination of angular velocity, $v_{\phi}(\rho)$}

In order to find the minimizer of the functional $S(\mathbf{v})$ of (\ref{eq:00}), we proceed as follows. First, for any distribution, $P(\rho)$, we will find the rotation-invariant velocity field, $\mathbf{v}_{P}$, which minimizes the functional
\be
S_{P}\left( \mathbf{v}\right)=\int_{R}^{\infty}\int_{0}^{2\pi}v_{\phi}^{2}(\rho)P(\rho)\,\rho d\rho\,d\phi
\label{eq:06a}
\ee
Secondly, we determine the distribution $P$ which minimizes the functional $S_P\left(\mathbf{v} \right)$. We consider the functional $S_{P}\left(\mathbf{\ast}\right)$ in the space of velocities with the given average tangential speed, $V$, that is
\be
\int_{R}^{\infty}\int_{0}^{2\pi}\frac{v_{\phi}(\rho)}{\rho}P(\rho)\,\rho
d\rho\,d\phi=\frac{V}{R}
\label{eq:06b}
\ee
Collecting (\ref{eq:06a}) and (\ref{eq:06b}) and introducing the Lagrange multiplier, $\lambda$, we can rewrite the auxiliary minimization problem in the following form:
\be
S_{P}^{\prime}=2\pi\int_{R}^{\infty}\big(v_{\phi}^{2}(\rho)\rho+\lambda v_{\phi}(\rho) \big)P(\rho)d\rho\rightarrow\mathrm{min}
\label{eq:06c}
\ee
Solution of the equation $\delta S_{P}^{\prime}=0$ gives
\be
v_{\phi}(\rho)=-\frac{\lambda}{\rho}
\label{eq:06d}
\ee
Substituting (\ref{eq:06d}) back into (\ref{eq:06b}), we get the expression for the Lagrange multiplier $\lambda$:
\be
\lambda=-\frac{V}{R}\left(2\pi\int_{R}^{\infty}\frac{P(\rho)}{\rho}
d\rho\right)^{-1}
\label{eq:06e}
\ee
Thus, we find the angular velocity, $v_{\phi}\equiv\left(\mathbf{v}_{P}\right)_{\phi}$, as a functional of $P(\rho)$:
\be
v_{\phi}=\frac{V}{2\pi R\rho}\left(\int_{R}^{\infty}\frac{P(\rho)}{\rho} d\rho\right)^{-1}
\label{eq:06g}
\ee

\section{Minimization of the action}

\subsection{General formalism}

To find the minimizer $P(\rho)$ of the action $S$, defined in (\ref{eq:00}), we substitute (\ref{eq:05}) and (\ref{eq:06g}) into (\ref{eq:00}), getting the following expression for $S\{P(\rho)\}$:
\begin{multline}
S\{P(\rho)\}=\int_{R}^{\infty}\big[v_{\rho}^{2}(\rho)+v_{\phi}(\rho)^{2}\big]P(\rho)\,\rho d\rho\,d\phi =2\pi\int_{R}^{\infty}\frac{D^{2}} {P(\rho)}\left(\frac{\partial P(\rho)}{\partial\rho}\right)^{2}\rho d\rho+ \\ 2\pi\frac{V^{2}}{4\pi^{2}R^{2}} \int_{R}^{\infty}
\frac{P(\rho)}{\rho}\left(\int_{R}^{\infty}\frac{P(\rho)}{\rho}d\rho\right)^{-2}d\rho =2\pi D^{2}\int_{R}^{\infty}\frac{1}{P(\rho)}\left(\frac{\partial P(\rho)}{\partial\rho}\right)^{2}\rho d\rho+\frac{V^{2}}{2\pi R^{2}}\left(\int_{R}^{\infty}\frac{P(\rho)}{\rho}d\rho\right)^{-1}
\label{eq:07}
\end{multline}
The equation for the normalization condition reads
\be
\int_{R}^{\infty}\int_{0}^{2\pi}P(\rho)\,\rho d\rho\,d\phi=1, 
\label{eq:07b}
\ee
Eq. \eq{eq:07b} should be added to the action $S\{P(\rho)\}$ with the Lagrange multiplier $\gamma$.

Let us make the substitution
\be
P(\rho)=Q^{2}(\rho)
\label{eq:07c}
\ee
and plug \eq{eq:07c} into (\ref{eq:07})--(\ref{eq:07b}). We obtain the functional
\be
S^{\prime}\{Q(\rho)\}=S\{Q^{2}(\rho)\}+2\pi\gamma\int_{R}^{\infty}Q^{2} (\rho)\,\rho d\rho. \label{eq:08}
\ee
Collecting all terms together, we can explicitly write (\ref{eq:08}) in the following form:
\be
S^{\prime}\{Q(\rho)\}=8\pi D^{2}\int_{R}^{\infty}\dot{Q}^{2}\rho d\rho
+\frac{V^{2}}{2\pi R^{2}}\left(  \int_{R}^{\infty}\frac{Q^{2}(\rho)}{\rho
}d\rho\right) ^{-1}+2\pi\gamma\int_{R}^{\infty}Q^{2}(\rho)\,\rho d\rho.
\label{eq:09}
\ee

In Eq.(\ref{eq:09}) we have used the notation $\dot{Q}=\frac{\partial Q(\rho)}{\partial\rho}$. It should be pointed out that the action (\ref{eq:09}) has the atypical second term, which, however, still allows us to proceed with the standard Euler-Lagrange minimization of the action (\ref{eq:09}). Equation $\delta S^{\prime}\{Q(\rho)\}=0$ leads to the following ordinary differential equation
\be
\frac{d}{d\rho}\left(\rho\frac{dQ(\rho)}{d\rho}\right)+\left(\frac{V^{2}}{16\pi^{2}D^{2} C^{2}R^{2}\rho} -\frac{\gamma\rho}{4D^{2}}\right) Q(\rho)=0
\label{eq:10}
\ee
where
\be
C=\displaystyle\int_{R}^{\infty}\frac{P(\rho)}{\rho}d\rho
\label{eq:10a}
\ee
Developing the first term in (\ref{eq:10}), we arrive at the boundary problem
\be
\begin{cases}
\displaystyle\frac{d^{2}Q(\rho)}{d\rho^{2}}+\frac{1}{\rho}\frac{dQ(\rho)}{d\rho}+
\frac{1}{4D^{2}}\left(  \frac{V^{2}}{4\pi^{2}C^{2}R^{2}\rho^{2}}-\gamma\right)
Q(\rho)=0 \medskip \\ Q(\rho=R)=0\medskip\\ Q(\rho\rightarrow\infty)\rightarrow 0
\end{cases}
\label{eq:11}
\ee

The constants $C$ and $\gamma$ should be determined self-consistently. We find $Q(\rho|C,\gamma)$ by solving the boundary problem (\ref{eq:11}), then we plug $P\equiv Q^{2}(\rho|C,\gamma)$ defined in (\ref{eq:07c}) into expression for $C$ getting
\be
C=\int_{R}^{\infty}\frac{Q^{2}(\rho|C,\gamma)}{\rho}d\rho\label{eq:12b}%
\ee
and $\gamma$ we determine from the normalization $\int_{R}^{\infty} \rho d\rho\int_{0}^{2\pi} d\phi\, P(\rho)=1$:
\be
2\pi\int_{R}^{\infty} Q^{2}(\rho|C,\gamma)\rho d\rho=1\label{eq:12c}%
\ee

\subsection{Analysis of the stationary asymptotic distribution}

We are interested in the asymptotic behavior of the stationary measure $P(\rho)$ in the vicinity of large circle of radius $R$, i.e. when $\rho=R+r$ and $0<r\ll R$. First, it is convenient to make in
(\ref{eq:11}) the substitution $Q(\rho)=U(\rho)\rho^{-1/2}$. In terms of the
function $U(\rho)$, equation (\ref{eq:11}) reads
\be
\frac{d^{2}U(\rho)}{d\rho^{2}}-\frac{1}{4}\left(\frac{\gamma}{D^{2}}-\frac{\Omega^{2}}
{\rho^{2}}\right) U(\rho)=0
\label{eq:12}
\ee
where $\Omega^{2}=\displaystyle1+\frac{V^{2}}{4\pi^{2}C^{2}D^{2}R^{2}}$. Writing $\rho=R+r$ in (\ref{eq:12}) and expanding this expression near $R$ up to linear terms in $r$, we get for $0\leq r\ll R$:
\be
\left\{\begin{array}[c]{l}
\displaystyle\frac{d^{2}U(r)}{dr^{2}}-\frac{1}{4}\left(\left(\frac{\gamma
}{D^{2}}-\frac{\Omega^{2}}{R^{2}}\right)+\frac{2\Omega^{2}}{R^{3}}r\right)
U(r)=0\medskip \\ U(r=0)=0 \medskip \\ U(r\rightarrow\infty)\rightarrow 0
\end{array} \right.
\label{eq:13}
\ee

Making use of the linear transform $r=u+vx$, we rewrite the boundary problem (\ref{eq:13}) as follows
\be
\begin{cases}
\displaystyle\frac{d^{2}U(x)}{dx^{2}}-xU(x)=0 \medskip \\
\displaystyle U\left(x=-\frac{u}{v}\right)=0 \medskip \\
\displaystyle U(x\rightarrow\infty)\rightarrow 0
\end{cases}
\label{eq:14}
\ee
where
\be
u=\frac{R}{2}\left(1-\frac{\gamma R^{2}}{D^{2}\Omega^{2}}\right); \qquad v=\frac{2^{1/3}R} {\Omega^{2/3}}
\label{eq:14a}
\ee
The general solution of (\ref{eq:14})--(\ref{eq:14a}) is
\be
U(x) = C_{1}\, \mathrm{Ai}\left(x+a_{1}+\frac{u}{v}\right)
\label{eq:15a}
\ee
where $\displaystyle \mathrm{Ai}(x) = \frac{1}{\pi}\int_{0}^{\infty} \cos\left(\frac{t^{3}}{3}+ tx\right)dt$. Hence,
\be
U(r) = C_{1}\, \mathrm{Ai}\left(  \frac{r-u}{v} + a_{1} + \frac{u}{v} \right)
\equiv C_{1}\, \mathrm{Ai}\left(  \frac{r}{v}+a_{1}\right)
\label{eq:15b}
\ee
where $C_{1}>0$ is the constant to be determined via (\ref{eq:12b}) and $a_{1}\approx-2.33811$ is the first (smallest in absolute value) zero of the Airy function, i.e. the solution of the equation $\mathrm{Ai}(a_{1})=0$. Correspondingly
\be
U(r) = C_{1}\, \mathrm{Ai}\left(  \left(  \frac{\Omega^{2}}{2}\right)^{1/3}\frac{r}{R}+
a_{1}\right)  \equiv C_{1}\, \mathrm{Ai}\left(\frac{\mu^{1/3}r}{R}+a_{1}\right)  \label{eq:16}%
\ee
where we have defined
\be
\mu= \frac{\Omega^{2}}{2}= \frac{1}{2}\left(1+\left(\frac{V}{2\pi C D R}\right)^{2}\right)  \label{eq:16a}
\ee

Using (\ref{eq:16}) and taking into account that $Q(r)=U(r)/\sqrt{R+r}$ (where $0\leq r\ll R$), we get the explicit expression for the requested distribution function, $Q(r)$:
\be
Q(r)=\frac{C_{1}}{\sqrt{R+r}}\,\mathrm{Ai}\left(  \frac{\mu^{1/3}\,r}{R}%
+a_{1}\right)  \label{eq:17}%
\ee
The constants $C_{1}, C$ are fixed by two auxiliary conditions (\ref{eq:12b}) and (\ref{eq:12c}), namely:
\begin{subequations}
\begin{align}
\displaystyle C  &  = C_{1}^{2}\int_{0}^{\infty}\frac{dr}{(R+r)^{2} }\,\mathrm{Ai}^{2}
\left(\frac{\mu^{1/3}\,r}{R}+a_{1}\right) \medskip \label{eq:18a} \\
1 &  = 2\pi C_{1}^{2}\int_{0}^{\infty}dr\,\mathrm{Ai}^{2}\left(\frac {\mu^{1/3}\,r}{R}+a_{1}\right)  \label{eq:18b}
\end{align}
\end{subequations}
Evaluating the integral in (\ref{eq:18b}), we determine the constant $C_{1}$:
\be
C_{1}^{2}=\frac{\mu^{1/3}}{2\pi\left(  \mathrm{Ai}^{\prime}(a_{1})\right)
^{2}R};\qquad\mathrm{Ai}^{\prime}(a_{1})=\frac{d\mathrm{Ai}(x)}{dx}
\bigg|_{x=a_{1}}
\label{eq:19}
\ee
The constant $C$ we find perturbatively expanding the denominator in (\ref{eq:18a}) in the power series in $\frac{r}{R}\ll1$. In the zero's order approximation we have
\be
\frac{1}{(R+r)^{2}}\approx\frac{1}{R^{2}},
\label{eq:expan}
\ee
which allows us to evaluate $C$ in (\ref{eq:18a}) in the zero's leading term in $0\ge\frac{r}{R}\ll1$:
\be
C\approx\frac{1}{2\pi R^{2}}
\label{eq:20}
\ee
The first-order correction to $C$ is derived in Appendix.

Collecting (\ref{eq:17}), (\ref{eq:19}) and (\ref{eq:20}), we get the final expression for the distribution function $P(\rho)$ where $\rho = R+r$:
\be
P(\rho)=\frac{\mu^{1/3}}{2\pi\left(\mathrm{Ai}^{\prime}(a_{1})\right)
^{2}R\,\rho}\, \mathrm{Ai}^{2}\left(\frac{\mu^{1/3}\,(\rho-R)}{R}+a_{1}\right)
\qquad(\rho\ge R),
\label{eq:21}
\ee
In \eq{eq:21} $\mu$ is obtained by substituting (\ref{eq:20}) into (\ref{eq:16a}):
\be
\mu=\frac{1}{2}\left(1+\frac{V^{2}R^{2}}{D^{2}}\right) \bigg|_{R\gg1}\approx \frac{V^{2}R^{2}}{2D^{2}}.
\label{eq:21a}
\ee

\section{Conclusion}

We have shown in the paper that the stationary radial distribution, $P(\rho)$ of the random walk which winds with the tangential velocity $V$ around the impenetrable disc of radius $R$ for $R\gg 1$ converges to the stationary distribution involving Airy function, given by an explicit expression
\be
P(\rho)=\frac{V^{2/3}}{2^{4/3}\pi\left(\mathrm{Ai}^{\prime}(a_{1})\right)^{2} D^{2/3} R^{1/3}\,\rho}\, \mathrm{Ai}^{2}\left(\left(\frac{V^2}{2D^2}\right)^{1/3} \frac{\rho-R}{R^{1/3}}+a_{1}\right)
\qquad(\rho\ge R),
\label{eq:23}
\ee
where $a_{1}\approx-2.33811$ is the first zero of the Airy function. Typical trajectories are localized in the circular strip $[R, R+ \delta R^{1/3}]$, where $\delta$ is the constant which depends on the parameters $D$ and $V$ and is independent on $R$.

It should be emphasized, that the presence of the impenetrable disc which restricts Brownian flights is crucial for the asymptotic behavior (\ref{eq:23}) and for the localization of trajectories within the strip of width $R^{1/3}$. It seems instructive to compare two models which look pretty similar. The first model represents the "stretched" random walk of $N=cR$ steps evading a disk of radius $R$ and in main respects repeats the problem considered in our paper (see also \cite{NPSVV} (definitely, on has $c>4$ for the lattice version on the square lattice). The second model discussed in \cite{HP}, studies the deviations from the typical "Wulf shape" of strongly inflated random loop of $N$ steps enclosing the algebraic area $A=cN^{2}$. At $N\gg 1$ the fluctuations in the first model scale as $N^{1/3}$, while remain Gaussian with the typical width $N^{1/2}$ in the second model. In both models the trajectories are pushed to very improbable tiny regions of the phase space, however the presence of a large deviation regime seems not to be a sufficient condition to affect the statistics and the convexity of the solid boundary on which trajectories recline, is also very important.

\begin{acknowledgments}

The authors are grateful for encouraging discussions with B. Meerson, K. Polovnikov and A. Valov. The research of S.S. was partially supported by the RSF grant No. 20-41-09009 and S.N. acknowledges the support of the Basis Foundation fellowship No. 19-1-1-48-1.

\end{acknowledgments}

\begin{appendix}

\section{}

Keeping the first-order corrections in \eq{eq:expan}, i.e. taking $\frac{1}{(R+r)^2} \approx \frac{1}{R^2}\left(1-\frac{2r}{R}\right)$, the selfconsistency equation \eq{eq:18a} becomes as follows (compare to \eq{eq:20})
\be
C\approx \frac{C_1^2}{R^2} \int_0^{\infty} dr\, \mathrm{Ai}^2\left(\frac{\mu^{1/3}\, r}{R}+a_1\right)-\frac{2C_1^2}{R^3} \int_0^{\infty} r dr\, \mathrm{Ai}^2\left(\frac{\mu^{1/3}\, r}{R}+a_1\right) =\frac{1}{2\pi R^2}\left(1 - \frac{4|a_1|}{3\pi \mu^{1/3}}\right)
\label{eq:20a}
\ee
Substituting in \eq{eq:20a} the general expression for $\mu$ extracted from \eq{eq:16a}, we get the following equation for the constant $C$
\be
2\pi C R^2 \approx 1 - \frac{2^{7/3}}{3}\frac{|a_1|}{\left(1+\disp \frac{V^2}{4\pi^2 C^2 D^2 R^2}\right)^{1/3}}
\label{eq:24}
\ee

By definition (see \eq{eq:12b}), the constant $C$ should be positive. Seeking the solution of \eq{eq:24} in the perturbative form, we may write $C$ as follows at $R\gg 1$:
\be
C \approx \frac{1}{2\pi R^2} + \frac{\kappa_1}{R^3} + \frac{\kappa_2}{R^4}
\label{eq:25}
\ee
where $\kappa_1$ and $\kappa_2$ are yet unknown coefficients which should be determined by solving \eq{eq:24} in the limit $R\gg 1$. Substituting the ansatz \eq{eq:25} into \eq{eq:24} we immediately find that $\kappa_1=0$ and $\kappa_2$ is defined by the expression
\be
\kappa_2 = -\frac{2^{4/3}}{3\pi} \frac{|a_1|D^2}{V^2}
\label{eq:26}
\ee
Thus, in the first-order expansion in $r/R\ll 1$, the constant $C$, which enters in the integral \eq{eq:12b}, reads
\be
C \approx \frac{1}{2\pi R^2} - \frac{2^{4/3}}{3\pi} \frac{|a_1|D^2}{V^2R^4}
\label{eq:27}
\ee

\end{appendix}

\end{document}